\documentclass{article}

\usepackage{amsmath}
\usepackage{amssymb}
\usepackage{amsthm}
\usepackage{amscd}

\newtheorem{thm}{Theorem}[section]
\newtheorem{prop}[thm]{Proposition}
\newtheorem{lem}[thm]{Lemma}
\newtheorem{cor}[thm]{Corollary}
\newtheorem{ithm}{Theorem}

\theoremstyle{definition}
\newtheorem{dfn}[thm]{Definition}

\theoremstyle{remark}
\newtheorem*{acknowledgments}{Acknowledgments}

\newcommand{\R}{\mathbb{R}}
\newcommand{\Z}{\mathbb{Z}}

\newcommand{\Ext}{\mathrm{Ext}}
\newcommand{\Hom}{\mathrm{Hom}}
\newcommand{\Tor}{\mathrm{Tor}}
\newcommand{\Sq}{\mathrm{Sq}}

\newcommand{\BB}{\mathcal{B}_{\mathrm{RR}} \times \mathcal{B}_{\mathrm{NS}}}
\newcommand{\I}{I_\mathrm{CS}}

\newcommand{\im}{\sqrt{-1}}
\renewcommand{\d}{\partial}
\renewcommand{\l}{\ell}
\def\h#1{ \widehat{#1} }

\def\abs#1{ \lvert{#1}\rvert }

\title{Differential characters and the Steenrod squares}

\author{Kiyonori Gomi}

\date{}

\begin{document}

\maketitle

\begin{abstract}
The groups of differential characters of Cheeger and Simons admit a natural multiplicative structure. The map given by the squares of degree $2k$ differential characters reduces to a homomorphism of ordinary cohomology groups. We prove that the homomorphism factors through the Steenrod squaring operation of degree $2k$. A simple application shows that five-dimensional Chern-Simons theory for pairs of $B$-fields is $SL(2, \Z)$-invariant on spin manifolds.
\end{abstract}


\section{Introduction}
\label{sec:introduction}

The group of \textit{Differential characters}, introduced by Cheeger and Simons \cite{Chee-S}, is a certain refinement of ordinary cohomology involving information of differential forms. From the beginning, differential characters enjoy numerous applications to geometry, topology and mathematical physics. 

\medskip

We recall the definition here: let $X$ be a smooth manifold. We denote the group of singular $p$-chains with coefficients in $\Z$ by $C_p(X) = C_p(X; \Z)$, and singular $p$-cycles by $Z_p(X) \subset C_p(X)$. A \textit{differential character}  of degree $\l$ is defined to be a homomorphism $\chi : Z_\l(X) \to \R/\Z$ such that there exists a differential $(\l+1)$-form $\omega$ satisfying $\chi(\d \tau) = \int_\tau \omega \mod \Z$ for all $\tau \in C_{\l+1}(X)$. The group of differential characters of degree $\l$ is denoted by $\h{H}^\l(X, \R/\Z)$.

It is known \cite{Chee-S} that there is a bilinear map on differential characters:
$$
\cup : \ \h{H}^{\l_1}(X, \R/\Z) \times \h{H}^{\l_2}(X, \R/\Z) \longrightarrow
\h{H}^{\l_1 + \l_2 + 1}(X, \R/\Z).
$$
The bilinear map $\cup$, which we call the \textit{cup product}, is associative. The cup product is also graded commutative in the sense that $\chi_1 \cup \chi_2 = (-1)^{(\l_1 + 1)(\l_2 + 1)} \chi_2 \cup \chi_1$ for $\chi_i \in \h{H}^{\l_i}(X, \R/\Z)$. 

Now we consider the quadratic map 
$$
\hat{q}^{\l} : \ \h{H}^\l(X, \R/\Z) \longrightarrow \h{H}^{2\l + 1}(X, \R/\Z)
$$ 
defined by $\hat{q}^{\l}(\chi) = \chi \cup \chi$. If we take an even integer $\l = 2k$, then $\hat{q}^{2k}$ gives rise to a homomorphism by the graded commutativity. Moreover, $\hat{q}^{2k}$ descends to give the following homomorphism of ordinary cohomology (see Lemma \ref{lem:q_hat_to_q}):
$$
q^{2k} : \ H^{2k+1}(X; \Z) \otimes \Z_2 \longrightarrow
\mathrm{Hom}(H_{4k+1}(X; \Z), \Z_2).
$$
In the case of $k = 0$ and $X = S^1$, an explicit formula of the cup product in \cite{Chee-S} leads to $q^0 = 1$, so that $q^{2k}$ is non-trivial in general.

The homomorphism $q^2$ has a relationship with the 5-dimensional topological field theory with Chern-Simons action studied by Witten in \cite{Wi}. To be more precise, we let $X$ be a compact oriented 5-dimensional manifold. We regard $\BB = \h{H}^2(X, \R/\Z) \times \h{H}^2(X, \R/\Z)$ as the space of fields (modulo gauge transformation), and consider the action functional $\I : \BB \to \R/\Z$ given by $\I (B_\mathrm{RR}, B_\mathrm{NS}) = - (B_\mathrm{RR} \cup B_\mathrm{NS})(X)$. Notice that $SL(2, \Z)$ acts on $\BB$ in the standard manner, however, $\I$ is not generally invariant under the action. This is the point $q^2$ appears: $\I$ has the $SL(2, \Z)$-invariance if and only if $q^2 = 0$. In the case where $X$ is the direct product of a 4-dimensional spin manifold and $S^1$, it is shown \cite{Wi} that $\I$ acquires the $SL(2, \Z)$-invariance on a certain subgroup in $\BB$.

\medskip

In the present paper, we prove generally that the homomorphism $q^{2k}$ factors through the \textit{Steenrod squaring operation} (\cite{St}, see \cite{Mil-S,Sp}): 
$$
\Sq^{2k} : \ H^{2k+1}(X; \Z_2) \longrightarrow H^{4k+1}(X; \Z_2).
$$ 
Let $\iota$ and $\pi$ be the homomorphisms in the universal coefficients theorems \cite{Sp}:
\begin{gather*}
0 \to
H^{2k+1}(X; \Z) \otimes \Z_2 \overset{\iota}{\to}
H^{2k+1}(X; \Z_2) \to
\Tor(H^{2k+2}(X; \Z), \Z_2) \to 0, \\
0 \to
\Ext(H_{4k}(X; \Z), \Z_2) \to
H^{4k+1}(X; \Z_2) \overset{\pi}{\to}
\Hom(H_{4k+1}(X; \Z), \Z_2) \to 0.
\end{gather*}
The main result of this paper is:
\begin{ithm} \label{ithm:main}
$q^{2k} = \pi \circ \Sq^{2k} \circ \iota$.
\end{ithm}

As a consequence, we can immediately see that the Chern-Simons action $\I$ has the $SL(2, \Z)$-invariance for a compact 5-dimensional spin manifold $X$. 

\medskip

As the group of differential characters refines ordinary cohomology theory, the notion of \textit{generalized differential cohomology} \cite{H-S} refines generalized cohomology theory involving information of differential forms. For differential (twisted) $K$-theory, a quantity similar to $q^{2k}$ is shown to be generally non-trivial \cite{F-M-S}. It may be interesting to establish counterparts of Theorem \ref{ithm:main} in generalized differential cohomology theories, which should be answered in a future work. 

\bigskip

The present paper is organized as follows. In Section \ref{sec:cup_prod_ordinary}, we review the cup product on ordinary cohomology and the Steenrod squaring operations. We also review homotopies between the cup product on cohomology with coefficients in $\R$ and the wedge product of differential forms. Such a homotopy is used in defining the cup product of differential characters. In Section \ref{sec:diff_char}, we describe the group of differential characters as a cohomology of a cochain complex, and introduce the cup product. After a study of the graded commutativity of the cup product, we prove Theorem \ref{ithm:main}.


\section{The cup product on ordinary cohomology}
\label{sec:cup_prod_ordinary}


\subsection{The cup product}

Let $X$ be a topological space. We denote by $C_p(X) = C_p(X; \Z)$ the group of singular $p$-chains on $X$, and by $\d : C_p(X) \to C_{p-1}(X)$ the boundary operator. The cycles and boundaries are denoted by $Z_p(X)$ and $B_p(X)$ as usual. 

We write $C_*(X)$ for the singular chain complex, and $C_*(X) \otimes_{\Z} C_*(X)$ for the chain complex obtained by the tensor product. Both chain complexes are augmented over the $\Z$-module $\Z$. We introduce a chain map 
$$
T : \ C_*(X) \otimes_{\Z} C_*(X) \to C_*(X) \otimes_{\Z} C_*(X)
$$ 
by $T(\sigma_1 \otimes \sigma_2) = (-1)^{\abs{\sigma_1}\abs{\sigma_2}}\sigma_2 \otimes \sigma_1$. 

By the \textit{method of acyclic models} (\cite{E-M}, see \cite{Fo, Sp}), we have: 

\begin{lem} \label{lem:acyclic_model_chain}
There exists a sequence $\{ D_i \}_{i \ge 0}$ of functorial homomorphisms $D_i : C_*(X) \to C_*(X) \otimes_\Z C_*(X)$ raising the degree by $i$ such that:
\begin{list}{}{\parsep=-2pt\topsep=4pt}
\item[(a)] 
$D_0$ is a chain map preserving the augmentations;

\item[(b)]
$\d D_i - (-1)^i D_i \d = D_{i-1} + (-1)^i TD_{i-1}$ \ 
for $i \ge 1$.
\end{list}
If $\{ D_i \}$ and $\{ D'_i \}$ are as above, then there exists a sequence $\{ E_i \}_{i \ge 0}$ of functorial homomorphisms $E_i : C_*(X) \to C_*(X) \otimes_\Z C_*(X)$ raising the degree by $i$ such that:
\begin{list}{}{\parsep=-2pt\topsep=4pt}
\item[(c)] 
$E_0 = 0$;

\item[(b)]
$D'_i - D_i = E_i + (-1)^iTE_i + \d E_{i+1} + (-1)^i E_{i+1} \d$ \ 
for $i \ge 0$.
\end{list}
\end{lem}

\begin{proof}
We follow \cite{Sp} (Chapter 5, Section 9). Let $R = \Z[t]/(t^2 -1)$ be the group ring of $\Z_2 = \Z/2\Z$. We define a chain complex $F_* = (F_*, \d)$ over $R$ as follows. For $i \ge 0$, $F_i = R \langle e_i \rangle$ is the free $R$-module of rank 1 generated by $e_i$. For $i < 0$, we put $F_i = 0$. The boundary operator $\d : F_i \to F_{i-1}$ is $\d(e_i) = (1 + (-1)^it) e_{i-1}$. We make $\Z$ into an $R$-module by letting $t \in R$ act as the identity. Note that $F_*$ is augmented over the $R$-module $\Z$.
The $R$-module structure on $F_*$ makes $F_* \otimes_\Z C_*(X)$ into a chain complex over $R$. We also make $C_*(X) \otimes_\Z C_*(X)$ into a chain complex over $R$ by letting $t \in R$ act as $T$. Both of these chain complexes are augmented over the $R$-module $\Z$.
Now, by the method of acyclic models, there exists a functorial chain map $D : F_* \otimes_\Z C_*(X) \to C_*(X) \otimes_\Z C_*(X)$ preserving the augmentations, and such chain maps $D$ and $D'$ are naturally chain homotopic. Then $D$ gives the sequence $\{ D_i \}_{i \ge 0}$ of functorial homomorphisms $D_i : C_*(X) \to C_*(X) \otimes_\Z C_*(X)$ stated in the lemma by setting $D(e_i \otimes \sigma) = D_i(\sigma)$. Similarly, a natural chain homotopy $E$ between $D$ and $D'$ gives $\{ E_i \}_{i \ge 0}$ by setting $E(e_i \otimes \sigma) = E_{i+1}(\sigma)$.
\end{proof}

Let $\Lambda$ be either $\Z, \Z_2$ or $\R$. We denote by $(C^p(X; \Lambda), \delta)$ the singular cochain complex with coefficients in $\Lambda$. Let $\{ D_i \}$ be a sequence of functorial homomorphisms in Lemma \ref{lem:acyclic_model_chain}. Using the natural multiplicative structure on $\Lambda$, we define a homomorphism of $\Lambda$-modules
$$
\cup : \ C^p(X; \Lambda) \otimes_\Lambda C^p(X; \Lambda) 
\longrightarrow C^{p+q}(X; \Lambda)
$$
by $f \cup g = D_0^*(f \otimes g)$. Since $D_0 : C_*(X) \to C_*(X) \otimes_\Z C_*(X)$ gives a \textit{diagonal approximation} (\cite{Sp}), $\cup$ induces the cup product on $H^*(X; \Lambda)$.

\medskip

For later convenience, we let $\{ D^i \}_{i \ge 0}$ be the sequence of functorial homomorphisms $D^i : C^*(X; \Lambda) \otimes_\Lambda C^*(X; \Lambda) \to C^*(X; \Lambda)$ given by $D^i(f \otimes g) = D_i^*(f \otimes g)$. By definition, $D^i$ lowers the degree by $i$, and satisfies 
$$
D^i \delta - (-1)^i \delta D^i = D^{i-1} + (-1)^i D^{i-1}T,
$$
where $T : C^*(X; \Lambda) \otimes_\Lambda C^*(X; \Lambda) \to C^*(X; \Lambda) \otimes_\Lambda C^*(X; \Lambda)$ is the cochain map given by $T(f \otimes g) = (-1)^{\abs{f} \abs{g}} g \otimes f$.


\subsection{The Steenrod squaring operations}

Recall that the \textit{Steenrod squaring operations} \cite{Mil-S,Sp,St} are the additive cohomology operations
$$
\Sq^i : \ H^*(X; \Z_2) \longrightarrow H^{*+i}(X; \Z_2)
$$
characterized by the following axioms:
\begin{list}{}{\parsep=-2pt\topsep=4pt}
\item[(a)] 
$\Sq^0 = 1$,

\item[(b)]
$\Sq^p(c) = c \cup c$ for $c \in H^p(X; \Z_2)$,

\item[(c)]
$\Sq^i(c) = 0$ for $c \in H^p(X; \Z_2)$ with $i > p$, 

\item[(d)]
$\Sq^i(c \cup c') = \sum_{i = j + j'} \Sq^j(c) \cup \Sq^{j'}(c')$ for $c, c' \in H^*(X; \Z_2)$.
\end{list}

We follow \cite{Sp} to realize Steenrod squaring operations: let $\{ D_i \}_{i \ge 0}$ be as in Lemma \ref{lem:acyclic_model_chain}, and $\{ D^i \}_{i \ge 0}$ the associated sequence of functorial homomorphisms $D^i : C^p(X; \Z_2) \otimes_{\Z_2} C^q(X; \Z_2) \to C^{p+q-i}(X; \Z_2)$. For $i \ge 0$ we define homomorphisms $\Sq^i : C^p(X; \Z_2) \longrightarrow C^{p+i}(X; \Z_2)$ by
$$
\Sq^i(c) =
\left\{
\begin{array}{cl}
0 & i > p, \\
D^{p - i}(c \otimes c) & i \le p.
\end{array}
\right.
$$
These homomorphisms induce the operations $\Sq^i : H^p(X; \Z_2) \to H^{p+i}(X; \Z_2)$.


\subsection{The cup product and the wedge product}

Let $X$ be a smooth manifold. The integration on singular simplices gives a functorial cochain map from the de Rham complex $(\Omega^*(X), d)$ to $(C^*(X; \R), \delta)$. As is well-known, there exists a homotopy between the wedge product $\omega_1 \wedge \omega_2$ and the cup product $\omega_1 \cup \omega_2$. In other words, the diagram: 
$$
\begin{CD}
\Omega^*(X) \otimes_\R \Omega^*(X) @>{\wedge}>> \Omega^*(X) \\
@VVV @VVV \\
C^*(X; \R) \otimes_\R C^*(X; \R) @>>{\cup}> C^*(X; \R)
\end{CD}
$$
is commutative up to a homotopy. In this subsection, we review such a homotopy, since it will be used in defining the cup product on the groups of differential characters (Definition \ref{dfn:cup_product}).

\begin{lem} \label{lem:acyclic_model_cochain}
There exists a sequence $\{ B^i \}_{i \ge 0}$ of functorial homomorphisms $B^i : \Omega^*(X) \otimes_\R \Omega^*(X) \to C^*(X; \R)$ lowering the degree by $i$ such that:
\begin{list}{}{\parsep=-2pt\topsep=4pt}
\item[(a)] 
$B^0 = 0$;

\item[(b)]
$\alpha \wedge \beta - \alpha \cup \beta = B^1 d (\alpha \otimes \beta) + \delta B^1(\alpha \otimes \beta)$ \ for $\alpha, \beta \in \Omega^*(X)$;

\item[(c)]
$- D^i = B^i + (-1)^i B^iT + B^{i+1}d + (-1)^i \delta B^{i+1}$ \ for $i \ge 1$.
\end{list}
\end{lem}

\begin{proof}
The proof is almost the same as that of Lemma \ref{lem:acyclic_model_chain}: we put $R' = \R[t]/(t^2 -1)$, and regard $\R$ as an $R'$-module by letting $t \in R'$ act as the identity. We define a cochain complex $(F^*, \delta)$ over $R'$ as follows. For $i \ge 0$, $F^i = R' \langle e^i \rangle$ is the free $R'$-module of rank 1. For $i < 0$, we put $F^i = 0$. The coboundary operator $\delta : F^i \to F^{i+1}$ is $\delta(e^i) = (1 - (-1)^it) e^{i+1}$. Then we have two cochain complexes $F^* \otimes_\R C^*(X; \R)$ and $\Omega^*(X) \otimes_\R \Omega_*(X)$ over $R'$. We define functorial cochain maps $W, U : \ \Omega^*(X) \otimes_\R \Omega^*(X) \to F^* \otimes_\R C^*(X; \R)$ by
\begin{align*}
W(\alpha \otimes \beta) &= (1+t)e^0 \otimes (\alpha \wedge \beta), \\
U(\alpha \otimes \beta) &= \sum_{i \ge 0} 
\left( e^i \otimes D^i(\alpha \otimes \beta) 
+ t e^i \otimes D^iT(\alpha \otimes \beta) \right).
\end{align*}
We now appeal to the method of acyclic models. (In particular, Theorem 7B in \cite{Fo} suits the present case.) Then there exists a natural cochain homotopy $B$ between $W$ and $U$. Because $B$ is an $R'$-module homomorphism, we have the following expression:
$$
B(\alpha \otimes \beta) = \sum_{i \ge 0}
\left( e^i \otimes B^{i+1}(\alpha \otimes \beta) 
+ t e^i \otimes B^{i+1}T(\alpha \otimes \beta) \right).
$$
We can easily verify that the sequence of homomorphisms $\{ B^i \}$ above has the properties stated in the present lemma.
\end{proof}


\section{Differential characters}
\label{sec:diff_char}


\subsection{Cohomology presentation}

We realize the group of differential characters as a cohomology group (\cite{D-F,H-S}).

\begin{dfn}
Let $X$ be a smooth manifold. For a positive integer $p$, we define $(\check{C}(p)^*, \check{d})$ to be the following cochain complex:
\begin{align*}
\check{C}(p)^q 
&= 
\left\{
\begin{array}{lc}
C^q(X; \Z) \times C^{q-1}(X; \R), & q < p, \\
C^q(X; \Z) \times C^{q-1}(X; \R) \times \Omega^q(X), & q \ge p, 
\end{array}
\right. \\
\check{d}(b, f) 
&= 
\left\{
\begin{array}{ll}
(\delta b, - b - \delta f), & \quad \quad
(b, f) \in \check{C}(p)^{q}, \ q < p-1, \\
(\delta b, - b - \delta f, 0), & \quad \quad
(b, f) \in \check{C}(p)^{p-1},
\end{array}
\right. \\
\check{d}(c, h, \omega) 
&= (\delta c, \omega - c - \delta h, d \omega), \quad \quad
(c, h, \omega) \in \check{C}(p)^q, \ q \ge p.
\end{align*}
We denote the cohomology group of this complex by $\check{H}(p)^q = \check{H}(p)^q(X)$.
\end{dfn}

The following lemma is easily shown.

\begin{lem} \label{lem:exact_sequences}
For a positive integer $p$, we have the exact sequences:
\begin{gather*}
0 \to
H^{p-1}(X; \R/\Z) \to
\check{H}(p)^p(X) \overset{\delta_1}{\to}
\Omega^{p}(X)_{\Z} \to 0, \\
0 \to
\Omega^{p-1}(X) / \Omega^{p-1}(X)_{\Z} \to
\check{H}(p)^p(X) \overset{\delta_2}{\to}
H^p(X; \Z) \to 0,
\end{gather*}
where $\Omega^q(X)_\Z$ is the group of closed integral q-forms.
\end{lem}

As in Section \ref{sec:introduction}, we denote by $\h{H}^\l(X, \R/\Z)$ the group of differential characters of degree $\l$.

\begin{lem}
There is an isomorphism $\check{H}(\l+1)^{\l+1}(X) \to \h{H}^\l(X, \R/\Z)$.
\end{lem}

\begin{proof}
For $x \in \check{H}(\l+1)^{\l+1}(X)$, let $(c, h, \omega) \in \check{Z}(\l+1)^{\l+1}$ be a representative of $x$. Note that $h \in C^{\l}(X; \R)$. We define a homomorphism $\chi : Z_\l(X) \to \R/\Z$ by $\chi(\sigma) = \langle h, \sigma \rangle \mod \Z$. Because $(c, h, \omega)$ is a cocycle, we have $\omega = c + \delta h$. Hence $\chi(\d \tau) = \int_\tau \omega \mod \Z$ for $\tau \in C_{\l+1}(X)$, so that $\chi \in \h{H}^\l(X, \R/\Z)$. The assignment $x \mapsto \chi$ gives rise to a well-defined homomorphism $\check{H}(\l+1)^{\l+1}(X) \to \h{H}^\l(X, \R/\Z)$. It is known \cite{Chee-S} that $\h{H}^\l(X, \R/\Z)$ fits into the exact sequences:
\begin{gather*}
0 \to
H^\l(X; \R/\Z) \to
\h{H}^\l(X, \R/\Z) \overset{\delta_1}{\to}
\Omega^{\l+1}(X)_{\Z} \to 0, \\
0 \to
\Omega^\l(X) / \Omega^\l(X)_{\Z} \to
\h{H}^\l(X, \R/\Z) \overset{\delta_2}{\to}
H^{\l+1}(X; \Z) \to 0.
\end{gather*}
Comparing the exact sequences above with those in Lemma \ref{lem:exact_sequences}, we can see that the homomorphism is bijective.
\end{proof}

In the remainder of this paper, we will mean by $\check{H}(\l+1)^{\l+1}(X)$ the group of differential characters $\h{H}^\l(X, \R/\Z)$.


\subsection{The cup product}
\label{subsec:cup_product_diff_char}

Now we introduce the cup product on $\check{H}(p)^p(X) \cong \h{H}^{p-1}(X, \R/\Z)$. By Lemma \ref{lem:acyclic_model_cochain}, we have a homotopy $B^1 : \Omega^p(X) \otimes_\R \Omega^q(X) \to C^{p+q-1}(X; \R)$ such that $\omega_1 \wedge \omega_2 - \omega_1 \cup \omega_2 = B^1d(\omega_1 \otimes \omega_2) + \delta B^1(\omega_1 \otimes \omega_2)$. 

\begin{dfn} \label{dfn:cup_product}
We define a homomorphism 
$$
\cup : \ \check{C}(p)^p \otimes_\Z \check{C}(q)^q \longrightarrow
\check{C}(p+q)^{p+q}
$$ 
by setting
\begin{multline*}
(c_1, h_1, \omega_1) \cup (c_2, h_2, \omega_2) \\
=
(c_1 \cup c_2, \
(-1)^p c_1 \cup h_2 + h_1 \cup \omega_2 + B^1(\omega_1 \otimes \omega_2), \
\omega_1 \wedge \omega_2).
\end{multline*}
The induced homomorphism $\cup : \check{H}(p)^p(X) \otimes_\Z \check{H}(q)^q(X) \to \check{H}(p+q)^{p+q}(X)$ is the cup product on differential characters.
\end{dfn}

It is known \cite{Chee-S} that the cup product is associative and graded commutative:
\begin{align*}
(x_1 \cup x_2) \cup x_3 &= x_1 \cup (x_2 \cup x_3), \\
x_1 \cup x_2 &= (-1)^{p_1 p_2} x_2 \cup x_1,
\end{align*}
where $x_i \in \check{H}(p_i)^{p_i}(X)$. The graded commutativity will be shown in the next subsection.

\medskip

It would be worth while to describe an example due to Cheeger and Simons \cite{Chee-S}. Notice that $\check{H}(1)^1(X) \cong \h{H}^0(X, \R/\Z) \cong C^\infty(X, U(1))$. We consider the case of $X = S^1$. For $f : S^1 \to U(1)$ we can find a map $F : \R \to \R$ such that $f(\theta) = \exp 2\pi\im F(\theta)$. We define $\Delta_f \in \Z$ by $F(\theta + 2\pi) = F(\theta) + \Delta_f$. Similarly, we introduce $G : \R \to \R$ to $g : S^1 \to U(1)$. Then the cup product $f \cup g \in \check{H}(2)^2(X) \cong \h{H}^1(S^1, \R/\Z)$ is expressed as
$$
(f \cup g) (S^1) = 
\Delta_f G(0) - \int_0^{2\pi} F \frac{dG}{d\theta} d\theta \mod \Z.
$$


\subsection{The graded commutativity of the cup product}

We here introduce an analogy of the sequence of homomorphisms $\{ D^i \}_{i \ge 0}$ to the cochain complex $(\check{C}(p)^*, \check{d})$. 

\begin{dfn}
Let $p, q, \bar{p}, \bar{q}$ be positive integers. 

(a) For a non-negative integer $i$, we define a homomorphism 
$$
F^i : \ 
\check{C}(p)^{\bar{p}} \otimes_\Z \check{C}(q)^{\bar{q}} \longrightarrow
C^{\bar{p} + \bar{q} - i - 1}(X; \R)
$$
as follows:
\begin{align*}
F^{2j}( (c_1, h_1, \omega_1) \otimes (c_2, h_2, \omega_2) ) 
&=
(-1)^{\bar{p}}D^{2j-1}(h_1 \otimes h_2) 
+ (-1)^{\bar{p}}D^{2j}(c_1 \otimes h_2) \\
&\quad
+ D^{2j}(h_1 \otimes \omega_2) 
+ B^{2j+1}(\omega_1 \otimes \omega_2), \\
F^{2j + 1}( (c_1, h_1, \omega_1) \otimes (c_2, h_2, \omega_2) ) 
&=
(-1)^{\bar{p}} D^{2j}(h_1 \otimes h_2) 
- D^{2j+1}(h_1 \otimes c_2) \\
&\quad
- (-1)^{\bar{p}} D^{2j+1}(\omega_1 \otimes h_2) 
- B^{2j+2}(\omega_1 \otimes \omega_2),
\end{align*}
where $D^{-1}(h_1 \otimes h_2) = 0$. When $\omega_1$ or $\omega_2$ is irrelevant, we substitute 0 for it.

(b) We define a sequence $\{ G^i \}_{i \ge 0}$ of homomorphisms
$$
G^i : \
\check{C}(p)^{\bar{p}} \otimes_\Z \check{C}(q)^{\bar{q}} \longrightarrow
\check{C}(p + q)^{\bar{p} + \bar{q} - i}
$$
by setting
\begin{multline*}
G^i((c_1, h_1, \omega_1) \otimes (c_2, h_2, \omega_2)) = \\
(D^i(c_1 \otimes c_2), \
F^i((c_1, h_1, \omega_1) \otimes (c_1, h_2, \omega_2)), \
W^i(\omega_1 \otimes \omega_2)),
\end{multline*}
where $W^i(\omega_1 \otimes \omega_2) = 0$ for $i > 0$ and $W^0(\omega_1 \otimes \omega_2) = \omega_1 \wedge \omega_2$. When $\omega_1$ or $\omega_2$ is irrelevant, we again substitute 0 for it.
\end{dfn}

\begin{prop} \label{prop:graded_commutativity}
The sequence of homomorphisms $\{ G^i \}_{i \ge 0}$ satisfies
\begin{list}{}{\parsep=-2pt\topsep=4pt}
\item[(a)] 
$G^0 \check{d} = \check{d} G^0$;

\item[(b)]
$G^i \check{d} - (-1)^i \check{d} G^i = G^{i-1} + (-1)^i G^{i-1}T$ \ 
for $i \ge 1$,
\end{list}
where we defined the cochain map $T : \ \check{C}(p)^{\bar{p}} \otimes_\Z \check{C}(q)^{\bar{q}} \to \check{C}(q)^{\bar{q}} \otimes_\Z \check{C}(p)^{\bar{p}}$ by setting $T(x_1 \otimes x_2) = (-1)^{\bar{p} \bar{q}}x_2 \otimes x_1$.
\end{prop}

\begin{proof}
We can directly prove (a) by Lemma \ref{lem:acyclic_model_cochain}. The proof of (b) amounts to showing the following formulae:
\begin{align*}
F^{2j} - F^{2j}T &= - D^{2j+1} - \delta F^{2j+1} + F^{2j+1}\check{d}, \\
F^{2j+1} + F^{2j+1}T & = D^{2j+2} + \delta F^{2j+2} + F^{2j+2} \check{d},
\end{align*}
where $j \ge 0$. We can show them by the help of Lemma \ref{lem:acyclic_model_chain} and Lemma \ref{lem:acyclic_model_cochain}.
\end{proof}

Since $G^0 : \check{C}(p)^p \otimes_\Z \check{C}(q)^q \to\check{C}(p + q)^{p + q}$ induces the cup product, we have:

\begin{cor}
The cup product $\cup : \check{H}(p)^p \otimes_\Z \check{H}(q)^q \to \check{H}(p+q)^{p+q}$ is graded commutative.
\end{cor}


\subsection{The main theorem}

Recall that we defined $\hat{q}^{\l} : \check{H}(\l+1)^{\l+1}(X) \longrightarrow \check{H}(2\l + 2)^{2\l+2}(X)$ by setting $\hat{q}^{\l}(x) = x \cup x$.

\begin{lem} \label{lem:q_hat_to_q}
For $k \ge 0$, the map $\hat{q}^{2k}$ induces the following homomorphism
$$
q^{2k} : \ H^{2k+1}(X; \Z) \otimes \Z_2 \longrightarrow
\mathrm{Hom}(H_{4k+1}(X; \Z), \Z_2).
$$
\end{lem}

\begin{proof}
Since the cup product is graded commutative, we have
$$
\hat{q}^{2k}(x + y) = x \cup x + x \cup y + y \cup x + y \cup y = x \cup x + y \cup y = \hat{q}^{2k}(x) + \hat{q}^{2k}(y).
$$
We focus on the exact sequences in Lemma \ref{lem:exact_sequences}. Clearly, we have $\delta_1(\hat{q}^{2k}(x)) = \delta_1(x) \wedge \delta_1(x) = 0$. For $\alpha \in \Omega^{2k}(X)/\Omega^{2k}(X)_\Z$, we also have $\hat{q}^{2k}(\alpha) = \frac{1}{2}d(\alpha \wedge \alpha) = 0$ in $\Omega^{4k+1}(X)/\Omega^{4k+1}(X)_\Z$. Hence the homomorphism $\hat{q}^{2k}$ descends to 
$$
\hat{q}^{2k} : \ H^{2k+1}(M; \Z) \longrightarrow
H^{4k+1}(M; \R/\Z) \cong \Hom(H_{4k+1}(M; \Z), \R/\Z).
$$
Using again the graded commutativity, we have $2 \hat{q}^{2k}(x) = 0$, so that $(\hat{q}^{2k}(x))(\sigma)$ belongs to $(\frac{1}{2}\Z)/\Z \subset \R/\Z$. The identification $(\frac{1}{2}\Z)/\Z \cong \Z/2\Z$ gives
$$
\hat{q}^{2k} : \ H^{2k+1}(M; \Z) \longrightarrow \Hom(H_{4k+1}(M; \Z), \Z_2).
$$
Because $(\hat{q}^{2k}(2x))(\sigma) = 2(\hat{q}^{2k}(x))(\sigma) = 0$, the $\hat{q}^{2k}$ descends to give $q^{2k}$.
\end{proof}

We are in the position of proving the main result:

\begin{thm} \label{thm:main}
$q^{2k} = \pi \circ \Sq^{2k} \circ \iota$.
\end{thm}

\begin{proof}
For any abelian group $A$, we have the natural monomorphism
$$
\Hom(H^{2k+1}(M; \Z) \otimes \Z_2, A) \longrightarrow 
\Hom(H^{2k+1}(M; \Z), A).
$$
So we will verify the coincidence of $q^{2k}$ and $\pi \circ \Sq^{2k} \circ \iota$ regarding them as homomorphisms $H^{2k+1}(M; \Z) \to \Hom(H_{4k+1}(X; \Z), \Z_2)$. Let $x = (c, h, \omega)$ be a cocycle in $\check{Z}(2k+1)^{2k+1}$. Because any cochain with coefficients in $\R$ is divisible by 2, the following formula is derived from Proposition \ref{prop:graded_commutativity}:
$$
x \cup x = 
\left(c \cup c, \
- \frac{1}{2} D^1(c \otimes c) 
- \delta \left(\frac{1}{2} F^1(x \otimes x) \right), \ 0 \right).
$$
Therefore we have the following expression of $q^{2k}$:
$$
\left( q^{2k} ([c]) \right) (\sigma)
= \langle D^1(c \otimes c), \sigma \rangle \mod 2\Z,
$$
where $\sigma \in Z_{4k+1}(X)$. We can see that this expression of $q^{2k}$ coincides with that of $\pi \circ \Sq^{2k} \circ \iota$, by means of the realization of $\Sq^{2k}$ in Section \ref{sec:cup_prod_ordinary}.
\end{proof}

An example of $q^{2k}$ is given by taking $X = S^1$. Then $\Sq^0 = 1$ implies that $q^0 : \Z_2 \to \Z_2$ is the identity map. We can verify this example directly by using the formula at the end of Subsection \ref{subsec:cup_product_diff_char}.

\medskip

In general, for a compact oriented $(4k+1)$-dimensional smooth manifold $X$ without boundary, the Steenrod squaring operation $\Sq^{2k}$ is expressed as $\Sq^{2k}(c) = v_{2k}(X) \cup c$ for $c \in H^{2k+1}(X; \Z_2)$, where $v_{2k}(X) \in H^{2k}(X; \Z_2)$ is the $2k$th \textit{Wu class} of $X$, (\cite{Mil-S,Sp}). Note that, in this case, the evaluation of the fundamental class $[X] \in H_{4k+1}(X; \Z)$ of $X$ simplifies $q^{2k}$ as 
$$
q^{2k} : \ H^{2k+1}(X; \Z) \otimes \Z_2 \longrightarrow \Z_2.
$$

\begin{cor}
For a compact oriented $(4k+1)$-dimensional smooth manifold $X$ without boundary, we can express $q^{2k} : H^{2k+1}(X, \Z) \otimes \Z_2 \to \Z_2$ as
$$
q^{2k}(c) = \langle v_{2k}(X) \cup \iota(c), [X] \rangle.
$$
\end{cor}

For example, when $X$ is 5-dimensional and spin, we have $v_{2}(X) = w_{2}(X) = 0$, so that $q^{2} = 0$.

\bigskip

\begin{acknowledgments}
The relation between $q^{2k}$ and $\Sq^{2k}$ was originally suggested by M. Furuta. I would like to express my gratitude to him. I would also like to thank T. Kohno, Y. Terashima, Y. Kametani and H. Sasahira for valuable discussions. The author's research is supported by Research Fellowship of the Japan Society for the Promotion of Science for Young Scientists.
\end{acknowledgments}


\begin{flushleft}
Graduate school of Mathematical Sciences, The University of Tokyo, \\
Komaba 3-8-1, Meguro-Ku, Tokyo, 153-8914 Japan. \\
e-mail: kgomi@ms.u-tokyo.ac.jp
\end{flushleft}

\end{document}